\documentclass{article}
\usepackage{amsmath}
\usepackage{amsmath,amssymb,latexsym,color}
\usepackage[mathscr]{eucal}
\newtheorem{thm}{Theorem}[section]

\newtheorem{lemma}[thm]{Lemma}
\newtheorem{cor}[thm]{Corollary}

\newtheorem{example}{Example}[section]
\newtheorem{defin}{Definition}[section]

\newtheorem{remark}{Remark}[section]

\newcommand{\proof}{{\it Proof.\quad}}
\newcommand{\qed}{\hfill\Box\medskip}

\usepackage{CJK}
\begin{document}

\renewcommand{\baselinestretch}{1.2}

\title{\bf Erd\H{o}s-Ko-Rado theorem and bilinear forms graphs for matrices
 over residue class rings}

\author{
Jun Guo\thanks{Corresponding author. guojun$_-$lf@163.com}\\
{\footnotesize Department of Mathematics, Langfang Normal University, Langfang  065000,  China} }
 \date{ }
 \maketitle

\begin{abstract}
Let $h=\prod_{i=1}^{t}p_i^{s_i}$ be
its decomposition into a product of powers of distinct primes, and $\mathbb{Z}_{h}$ be the residue class ring modulo $h$. Let $1\leq r\leq m\leq n$ and $\mathbb{Z}_{h}^{m\times n}$ be the set of all $m\times n$ matrices over $\mathbb{Z}_{h}$. The generalized bilinear forms graph over  $\mathbb{Z}_{h}$, denoted by $\hbox{Bil}_r(\mathbb{Z}_{h}^{m\times n})$, has the vertex set $\mathbb{Z}_{h}^{m\times n}$, and two distinct vertices $A$ and $B$ are adjacent if the inner rank of $A-B$ is less than or equal to $r$.  In this paper, we determine the clique number
and geometric structures of maximum cliques of $\hbox{Bil}_r(\mathbb{Z}_{h}^{m\times n})$.
As a result,  the Erd\H{o}s-Ko-Rado theorem for $\mathbb{Z}_h^{m\times n}$ is obtained.



\end{abstract}

\section{Introduction}
Let $\mathbb{Z}$ denote  the integer ring.
For $a,b,h\in \mathbb{Z}$, integers $a$ and $b$ are said to be {\it congruent modulo} $h$
 if $h$ divides $a-b$, and  denoted by $a\equiv b\; {\rm mod}\; h$.
  Suppose that $h=\prod_{i=1}^{t}p_i^{s_i}$ is its decomposition into a product of powers of distinct primes. Let $\mathbb{Z}_{h}$ denote the residue class ring modulo $h$ and $\mathbb{Z}_{h}^\ast$ denote its unit group. Then $\mathbb{Z}_{h}$ is a  principal ideal ring and $|\mathbb{Z}_{h}^\ast|=h\prod_{i=1}^t(1-p_i^{-1})$. Note that $(p_i)$, where $i=1,2,\ldots,t$, are all the maximal ideals of $\mathbb{Z}_{h}$.
By \cite{Ireland}, $\mathbb{Z}_{h}\cong\mathbb{Z}_{p_1^{s_1}}\oplus\mathbb{Z}_{p_2^{s_2}}\oplus\cdots\oplus\mathbb{Z}_{p_t^{s_t}}$ and
$\mathbb{Z}_{h}^\ast\cong\mathbb{Z}_{p_1^{s_1}}^\ast\times\mathbb{Z}_{p_2^{s_2}}^\ast\times\cdots\times\mathbb{Z}_{p_t^{s_t}}^\ast$.
For $a\in\mathbb{Z}$, we denote also by $a$ the congruence class of $a$ modulo $h$.

For a subset $S$ of $\mathbb{Z}_{h}$, let $S^{m\times n}$ be the set of all $m\times n$ matrices over $S$.
 Let ${}^t\!A$ denote the transpose matrix of a matrix $A$ and $\hbox{det}(X)$ the determinant of a square matrix $X$ over $\mathbb{Z}_{h}$. Let $I_r$ ($I$ for short) be the $r\times r$ identity matrix, and $0_{m,n}$ ($0$ for short) the $m\times n$ zero matrix.  Let  $\hbox{diag}(A_1,A_2,\ldots, A_k)$ denote the block diagonal matrix whose blocks along the main diagonal are matrices $A_1,A_2,\ldots, A_k$. The set of $n\times n$ invertible matrices forms a group under matrix multiplication, called the {\it general linear group} of degree $n$ over $\mathbb{Z}_{h}$ and denoted by $G\!L_n(\mathbb{Z}_{h})$. For $A\in\mathbb{Z}_{h}^{n\times n}$, by Corollary~2.21 in \cite{Brown},
$A\in G\!L_n(\mathbb{Z}_{h})$ if and only if $\hbox{det}(A)\in\mathbb{Z}_{h}^\ast.$

Let $A\in \mathbb{Z}_{h}^{m\times n}$ be a non-zero matrix. By Cohn's definition \cite{Cohn2}, the {\it inner rank} of $A$, denoted by $\rho(A)$, is the least integer $r$ such that
$A = BC$ where $B \in \mathbb{Z}_{h}^{m\times r}$ and $C \in \mathbb{Z}_{h}^{r\times n}$.
Let $\rho(0) =0$. For $A\in \mathbb{Z}_{h}^{m\times n}$, it is obvious that $\rho(A)\leq\min\{m,n\}$ and $\rho(A)=0$ if and only if $A=0$. For matrices over $\mathbb{Z}_{h}$, by \cite{Cohn,Cohn2},
$\rho(A)=\rho(SAT)$ where $S\in G\!L_{m}(\mathbb{Z}_{h})$ and $T\in G\!L_{n}(\mathbb{Z}_{h})$, $\rho(AB) \leq\min\{\rho(A),\rho(B)\}$ and
$$
\rho\left(\begin{array}{cc}
   A_{11} & A_{12}\\
   A_{21}  & A_{22}
   \end{array}\right) \geq\max\{\rho(A_{ij}) :1\leq i,j\leq 2\}.$$

The Erd\H{o}s-Ko-Rado theorem \cite{Erdos,Wilson}
 is a classical result in extremal set theory which obtained an upper bound on the size of a family of $m$-subsets of a set that every pairwise intersection has size at least $r$ and describes exactly which families meet this bound. The results on Erd\H{o}s-Ko-Rado theorem have inspired much research \cite{Frankl,Godsi2,Tanaka}. Let $1\leq r\leq m\leq n$. A family ${\cal F}\subseteq\mathbb{F}_q^{m\times n}$ is called $r$-{\it intersecting}
 if ${\rm rank}(A-B)\leq r$ for all $A,B\in{\cal F}$, where $\mathbb{F}_q^{m\times n}$ is the set of all $m\times n$
 matrices over the $q$-element finite field $\mathbb{F}_q$. Huang \cite{Huang-T} obtained an upper bound on the size of an $r$-intersecting family in $\mathbb{F}_q^{m\times n}$ and describes exactly which families meet this bound.
 As a natural extension, a family ${\cal F}\subseteq\mathbb{Z}_h^{m\times n}$ is called $r$-{\it intersecting} if $\rho(A-B)\leq r$ for all $A,B\in{\cal F}$.

The bilinear forms graph over a finite field plays an important role in geometry and combinatorics, and it has been extensively studied, see \cite{Brouwer,Wan3,Wang2}. As a natural extension, the {\it generalized bilinear forms graph} over  $\mathbb{Z}_{h}$, denoted by $\hbox{Bil}_r(\mathbb{Z}_{h}^{m\times n})$, has the vertex set $\mathbb{Z}_{h}^{m\times n}$, and two distinct vertices $A$ and $B$ are adjacent if   $\rho(A-B)\leq r$, where  $1\leq r\leq m\leq n$.  Note that $\hbox{Bil}_1(\mathbb{Z}_{h}^{m\times n})$ is the bilinear forms graph $\hbox{Bil}(\mathbb{Z}_{h}^{m\times n})$.
When $t=1$, Huang et al. \cite{Huang,Huang2} determined the clique number, the independence number and the chromatic number of $\hbox{Bil}_r(\mathbb{Z}_{h}^{m\times n})$.

Let $V(\Gamma)$ denote the vertex set of a graph $\Gamma$. For $A,B\in V(\Gamma)$, we write $A\sim B$ if vertices $A$ and $B$ are adjacent. A {\it clique} of a graph $\Gamma$ is a complete subgraph of $\Gamma$. A clique ${\cal C}$ is {\it maximal} if there is no clique of $\Gamma$ which properly contains ${\cal C}$ as a subset. A {\it maximum clique} of $\Gamma$ is a clique of $\Gamma$ which has maximum cardinality. The {\it clique number} $\omega(\Gamma)$ of $\Gamma$ is the number of vertices in a maximum clique.
An {\it independent set} of a graph $\Gamma$ is a subset of vertices such that no two vertices are adjacent. A {\it largest independent set} of $\Gamma$ is an independent set of maximum cardinality. The {\it independence number} $\alpha(\Gamma)$ is the number of vertices in a largest independent set of $\Gamma$.

An $\ell$-{\it coloring} of a graph $\Gamma$ is a homomorphism from $\Gamma$ to the complete graph $K_\ell$.
The {\it chromatic number} $\chi(\Gamma)$ of $\Gamma$ is the least value $k$ for which $\Gamma$ can be $k$-colored.
A graph $\Gamma$ is a {\it core} \cite{Godsil} if every endomorphism of $\Gamma$ is an automorphism. A
 subgraph $\Delta$ of a graph $\Gamma$ is a core of $\Gamma$ if it is a core and there exists some homomorphism from
 $\Gamma$ to $\Delta$. Every graph $\Gamma$ has a core, which is an induced subgraph and is unique up to
 isomorphism, see Lemma~6.2.2 in \cite{Godsil}.

The Smith normal forms of matrices over $\mathbb{Z}_{h}$ for $t=1,2$ are determined in \cite{Guo2,Newman}.
 In Section~2, we determine the Smith normal forms of matrices over $\mathbb{Z}_{h}$, and present some useful results for later reference.
 Note that every maximum clique of
 $\hbox{Bil}_r(\mathbb{Z}_{h}^{m\times n})$ is a largest $r$-intersecting family in $\mathbb{Z}_{h}^{m\times n}$ and vice versa.
  In Section~3, we determine the clique number, the independence number
  and the chromatic number of the generalized bilinear forms graph $\hbox{Bil}_r(\mathbb{Z}_{h}^{m\times n})$, and show that cores of  both $\hbox{Bil}_r(\mathbb{Z}_{h}^{m\times n})$  and its complement are maximum cliques.
  As a result,  the Erd\H{o}s-Ko-Rado theorem for $\mathbb{Z}_h^{m\times n}$ is obtained.

\section{Smith normal forms and ranks of matrices}
  Let $J_{(\alpha_1,\alpha_2,\ldots,\alpha_t)}=(\prod_{i=1}^{t}p_i^{\alpha_i})$, where $0\leq \alpha_i\leq s_i$ for $i=1,2,\ldots,t$. For brevity, we write $J_{(\alpha_1)}$ as $J_{\alpha_1}$ if $t=1$.

\begin{lemma}\label{lem2.1}
 Let $I$ be a nonempty proper subset  of $[t]=\{1,2,\ldots,t\}$. Suppose that $\alpha_i\geq s_i$ for each $i\in I$, and  $0\leq \alpha_j<s_j$ for each $j\in[t]\setminus I$. Then there exists some unit $u$  such that $u\prod_{\theta=1}^tp_\theta^{\alpha_\theta}=\prod_{i\in I}p_i^{s_i}\prod_{j\in[t]\setminus I}p_j^{\alpha_j}$.
\end{lemma}
\proof Without loss of generality, we may assume that $I=\{1,2,\ldots,|I|\}$. Note that $p_\xi+\prod_{i\in I\setminus\{\xi\}}p_i\prod_{j=|I|+1}^tp_j^{s_j-\alpha_j}$
 is a unit for each $\xi\in I$. From
 $$\left(p_1+\prod_{i=2}^{|I|}p_i\prod_{j=|I|+1}^{t}p_j^{s_j-\alpha_j}\right)
 p_1^{s_1}\prod_{\theta=2}^tp_\theta^{\alpha_\theta}
=p_1^{1+s_1}\prod_{\theta=2}^tp_\theta^{\alpha_\theta},$$
 we deduce that
there exists some unit $u_1$ such that $u_1\prod_{\theta=1}^tp_\theta^{\alpha_\theta}
=p_1^{s_1}\prod_{\theta=2}^tp_\theta^{\alpha_\theta}.$
Similarly, from
$$\left(p_2+\prod_{i\in I\setminus\{2\}}p_i\prod_{j=|I|+1}^{t}p_j^{s_j-\alpha_j}\right)
 p_1^{s_1}p_2^{s_2}\prod_{\theta=3}^tp_\theta^{\alpha_\theta}
=p_1^{s_1}p_2^{1+s_2}\prod_{\theta=3}^tp_\theta^{\alpha_\theta},$$
we deduce that there exists some unit $u_2$ such that $u_2 p_1^{s_1}\prod_{\theta=2}^tp_\theta^{\alpha_\theta}
=p_1^{s_1}p_2^{s_2}\prod_{\theta=3}^tp_\theta^{\alpha_\theta}.$
And so on, there exists some unit $u_{|I|}$ such that
$$u_{|I|} \prod_{i=1}^{|I|-1}p_i^{s_i}\prod_{\theta=|I|}^tp_\theta^{\alpha_\theta}
=\prod_{i=1}^{|I|}p_i^{s_i}\prod_{\theta=|I|+1}^tp_\theta^{\alpha_\theta}.$$
Let $u=\prod_{i=1}^{|I|}u_i$. Then $$u\prod_{\theta=1}^tp_\theta^{\alpha_\theta}=\prod_{i=1}^{|I|}p_i^{s_i}\prod_{\theta=|I|+1}^tp_\theta^{\alpha_\theta}.$$
Therefore, the  desired result follows. $\qed$

\begin{lemma}\label{lem2.2}
 The principal ideals $J_{(\alpha_1,\alpha_2,\ldots,\alpha_t)}$, where $0\leq \alpha_{i}\leq s_i$ for $i=1,2,\ldots,t$, are all the ideals of $\mathbb{Z}_{h}$.
\end{lemma}
\proof Since the ring $\mathbb{Z}_{h}$ is a principal ideal ring, each ideal has the form $(x)$. If $x=0$, then  $(0)=J_{(s_1,s_2,\ldots,s_t)}$. Suppose $1\leq x< \prod_{i=1}^tp_i^{s_i}$. Then there exist the unique $u\in\mathbb{Z}$ and the unique vector $(\beta_1,\beta_2,\ldots,\beta_t)$  such that $x=u\prod_{i=1}^tp_i^{\beta_i}$ and $u\in\mathbb{Z}_{h}^\ast$.
Since $x\not=0$ in $\mathbb{Z}_{h}$, there exists some $\beta_i$ such that $\beta_i<s_i$. By Lemma~\ref{lem2.1},
there exists some unit $v$ such that $v\prod_{i=1}^tp_i^{\beta_i}=\prod_{i=1}^tp_i^{\min\{\beta_i,s_i\}}$. Let $\alpha_i=\min\{\beta_i,s_i\}$
for $i=1,2,\ldots,t$. Then $0\leq\alpha_i\leq s_i$
for $i=1,2,\ldots,t$, and $(x)=(\prod_{i=1}^tp_i^{\beta_i})=(\prod_{i=1}^tp_i^{\alpha_i})$. $\qed$

\begin{lemma}\label{lem2.3}
Every non-zero element $x$ in $\mathbb{Z}_{h}$ can be written as $x=u\prod_{i=1}^tp_i^{\alpha_i}$, where $u$ is a unit, and
$0\leq\alpha_i\leq s_i$ for $i=1,2,\ldots,t$.
Moreover, the vector $(\alpha_1,\alpha_2,\ldots,\alpha_t)$ is unique
and $u$ is unique modulo the ideal $J_{(s_1-\alpha_1,s_2-\alpha_2,\ldots,s_t-\alpha_t)}$.
\end{lemma}
\proof Similar to the proof of Lemma~\ref{lem2.2}, there exist the unique $u_1\in\mathbb{Z}$ and the unique vector $(\beta_1,\beta_2,\ldots,\beta_t)$ such that  $x=u_1\prod_{i=1}^tp_i^{\beta_i}$ and $u_1\in\mathbb{Z}_h^\ast$. Since $x\not=0$, there exists some $i$ such that $\beta_i<s_i$. Write $\alpha_i=\min\{\beta_i,s_i\}$ for $i=1,2,\ldots,t$. By Lemma~\ref{lem2.1}, there exists some unit $u_2$ such that $$x=u_1\prod_{i=1}^tp_i^{\beta_i}=u_1u_2\prod_{i=1}^tp_i^{\alpha_i}.$$

Suppose $x=u\prod_{i=1}^tp_i^{\alpha_i}=v\prod_{i=1}^tp_i^{\gamma_i}$, where $u$ and $v$ are units, and
$0\leq\alpha_i,\gamma_i\leq s_i$ for $i=1,2,\ldots,t$. We claim that $\alpha_i=\gamma_i$ for $i=1,2,\ldots,t$.
If there exists some $i$ such that $\alpha_i\not=\gamma_i$, without loss of generality, we may assume that
$\alpha_i>\gamma_i$. Since $h$ divides $u\prod_{i=1}^tp_i^{\alpha_i}-v\prod_{i=1}^tp_i^{\gamma_i}$ and $\alpha_i\leq s_i$, we deduce that $p_i^{\alpha_i}$ divides $u\prod_{i=1}^tp_i^{\alpha_i}-v\prod_{i=1}^tp_i^{\gamma_i}$, which implies that
$p_i$ divides $v\prod_{j\not=i}p_j^{\gamma_j}$, a contradiction.

Suppose $x=u\prod_{i=1}^tp_i^{\alpha_i}=v\prod_{i=1}^tp_i^{\alpha_i}$,  where $u$ and $v$ are units. Then $h$ divides $(u-v)\prod_{i=1}^tp_i^{\alpha_i}$,
which implies that $\prod_{i=1}^tp_i^{s_i-\alpha_i}$ divides $u-v$. It follows that $u-v\in(\prod_{i=1}^tp_i^{s_i-\alpha_i})$.
Therefore, $u$ is unique modulo the ideal $J_{(s_1-\alpha_1,s_2-\alpha_2,\ldots,s_t-\alpha_t)}$.
$\qed$

Let $a,b\in \mathbb{Z}_{h}$. If there is an element $c\in \mathbb{Z}_{h}$ such that $b=ac$, we say that $a$ {\it divides} $b$ and denoted by $a|b$. Two elements $a$ and $b$ are said to be {\it associates} if $a=ub$ for some unit $u\in\mathbb{Z}_{h}^\ast$.

\begin{lemma}\label{lem2.4}
For $a,b\in\mathbb{Z}_h$,  $a$ and $b$ are associates if and only if $\pi_i(a)$ and $\pi_i(b)$ are associates
for $i=1,2,\ldots,t$.
\end{lemma}
\proof Suppose that $a$ and $b$ are associates. Then there exists some unit $u$ such that $a=bu$.
Therefore, we have $\pi_i(a)=\pi_i(bu)=\pi_i(b)\pi_i(u)$ for $i=1,2,\ldots,t$, which imply that
$\pi_i(a)$ and $\pi_i(b)$ are associates
for $i=1,2,\ldots,t$. Conversely, assume that $\pi_i(a)$ and $\pi_i(b)$ are associates
for $i=1,2,\ldots,t$. Then there exist $u_i\in\mathbb{Z}_{p_i^{s_i}}^\ast$ such that
$\pi_i(a)=\pi_i(b)u_i$ for $i=1,2,\ldots,t$. Since $\mathbb{Z}_{h}^\ast\cong\mathbb{Z}_{p_1^{s_1}}^\ast\times\mathbb{Z}_{p_2^{s_2}}^\ast
\times\cdots\times\mathbb{Z}_{p_t^{s_t}}^\ast$, there exists the unique $u\in\mathbb{Z}_{h}^\ast$
such that $\pi_i(u)=u_i$  for $i=1,2,\ldots,t$. Therefore, we have $\pi_i(a-bu)=0$ for $i=1,2,\ldots,t$, which imply that $a=bu$. Thus $a$ and $b$ are associates. $\qed$

 The matrices $A$ and $B$ in $\mathbb{Z}_{h}^{m\times n}$ are called {\it equivalent} if $S^{-1}AT=B$ for some $S\in G\!L_{m}(\mathbb{Z}_{h})$ and
 $T\in G\!L_{n}(\mathbb{Z}_{h})$. The direct product $G\!L_{m}(\mathbb{Z}_{h})\times G\!L_{n}(\mathbb{Z}_{h})$ acts on the set $\mathbb{Z}_{h}^{m\times n}$ in the following way:
$$\begin{array}{rcl}
\mathbb{Z}_{h}^{m\times n}\times(G\!L_{m}(\mathbb{Z}_{h})\times G\!L_{n}(\mathbb{Z}_{h})) &\rightarrow & \mathbb{Z}_{h}^{m\times n},\\
(A,(S,T)) &\mapsto&S^{-1}AT.
\end{array}$$
Clearly, both $A$ and $B$ belong to the same orbit of $\mathbb{Z}_{h}^{m\times n}$ under the action of $G\!L_{m}(\mathbb{Z}_{h})\times G\!L_{n}(\mathbb{Z}_{h})$ if and only if they are equivalent.

\begin{lemma}\label{lem3.2}{\rm(See \cite{Newman}.)}
Let $m\leq n$. Then every matrix $A\in\mathbb{Z}_{p^s}^{m\times n}$ is equivalent to
$$
{\rm diag}\,(p^{\alpha_{1}},p^{\alpha_{2}},\ldots,p^{\alpha_{m}}),
$$
where $0\leq\alpha_{1}\leq\alpha_{2}\leq\cdots\leq\alpha_{m}\leq s$.
Moreover, the vector
$(\alpha_{1},\alpha_{2},\ldots,\alpha_{m})$ is uniquely determined by $A$.
\end{lemma}

Let  $\pi_i$  be the natural surjective homomorphism from
$\mathbb{Z}_{h}$ to $\mathbb{Z}_{p_i^{s_i}}$ for $i=1,2,\ldots,t$. For $A=(a_{uv})\in\mathbb{Z}_{h}^{m\times n}$, let
$\pi_i(A)=(\pi_i(a_{uv}))$ for $i=1,2,\ldots,t$. Then there is a bijective map $\pi$ from
$\mathbb{Z}_{h}^{m\times n}$ to $\mathbb{Z}_{p_1^{s_1}}^{m\times n}\times \mathbb{Z}_{p_2^{s_2}}^{m\times n}
 \times\cdots\times\mathbb{Z}_{p_t^{s_t}}^{m\times n}$ such that $\pi(A)=(\pi_1(A),\pi_2(A),\ldots,\pi_t(A))$ for every $A\in \mathbb{Z}_{h}^{m\times n}$.

\begin{lemma}\label{lem3.1}{\rm(See \cite{Rosen}.)}
Let  $\pi_i$ (resp. $\pi$) denote also by the restriction of $\pi_i$ (resp. $\pi$) on $G\!L_{n}(\mathbb{Z}_{h})$  for $i=1,2,\ldots,t$.
Then $\pi_i$ is a natural surjective homomorphism  from
$G\!L_{n}(\mathbb{Z}_{h})$ to $G\!L_{n}(\mathbb{Z}_{p_i^{s_i}})$ for $i=1,2,\ldots,t$, and
$\pi$ is an isomorphism from $G\!L_{n}(\mathbb{Z}_{h})$ to $G\!L_{n}(\mathbb{Z}_{p_1^{s_1}})\times G\!L_{n}(\mathbb{Z}_{p_2^{s_2}})\times\cdots\times G\!L_{n}(\mathbb{Z}_{p_t^{s_t}})$.
\end{lemma}

\begin{lemma}\label{lem3.3}{\rm(See Theorem~15.24 in \cite{Brown}.)}
Let $m\leq n$. Then every matrix $A\in\mathbb{Z}_{h}^{m\times n}$ is equivalent to
${\rm diag}\,(d_1,d_2,\ldots,d_m),$
where $d_j|d_{j+1}$ for $j=1,2,\ldots,m-1$. Moreover, if $A$ is also equivalent to
${\rm diag}\,(b_1,b_2,\ldots,b_m),$ where $b_j|b_{j+1}$ for $j=1,2,\ldots,m-1$, then $d_c$ and $b_c$ are associates for $c=1,2,\ldots,m$.
\end{lemma}

\begin{thm}{\rm (Smith normal form.)}\label{lem3.4}
Let $m\leq n$. Then every matrix $A\in\mathbb{Z}_{h}^{m\times n}$ is equivalent to
\begin{equation}\label{equa2}
D={\rm diag}\,\left(\prod_{i=1}^tp_i^{\alpha_{i1}},\prod_{i=1}^tp_i^{\alpha_{i2}},\ldots,\prod_{i=1}^tp_i^{\alpha_{im}}\right),
\end{equation}
where $0\leq\alpha_{i1}\leq\alpha_{i2}\leq\cdots\leq\alpha_{im}\leq s_i$ for $i=1,2,\ldots,t$.
Moreover, the array
\begin{equation}\label{equa3}
\Omega=((\alpha_{11},\alpha_{12},\ldots,\alpha_{1m}),(\alpha_{21},\alpha_{22},\ldots,\alpha_{2m}),\ldots,(\alpha_{t1},\alpha_{t2},\ldots,\alpha_{tm}))
\end{equation} is uniquely determined by $A$.
\end{thm}
\proof
For each $i$ with $1\leq i\leq t$, by Lemma~\ref{lem3.2},  there exist invertible matrices $S_i\in G\!L_m(p_i^{s_i})$ and $T_i\in G\!L_n(p_i^{s_i})$ such that
$$S_i\pi_i(A)T_i={\rm diag}\,(p_i^{\alpha_{i1}},p_i^{\alpha_{i2}},\ldots,p_i^{\alpha_{im}})=D_i,$$
where $0\leq\alpha_{i1}\leq\alpha_{i2}\leq\cdots\leq\alpha_{im}\leq s_i$.

By Lemma~\ref{lem3.1}, there exists the unique pair of matrices $(S,T)\in G\!L_m(\mathbb{Z}_h)\times G\!L_n(\mathbb{Z}_h)$ such that
$\pi_{i}(S)=S_i$ and $\pi_{i}(T)=T_i$ for $i=1,2,\ldots,t$. Therefore, we have
$$\pi_{i}(S)\pi_i(A)\pi_{i}(T)-D_i=0\;\hbox{for}\; i=1,2,\ldots,t.$$
 Write $S=(s_{uv}),T=(t_{uv})$ and $A=(a_{uv})$.
For $i=1,2,\ldots,t$, from $\pi_{i}(S)\pi_i(A)\pi_{i}(T)-D_i=0$, we deduce that
$$\pi_i\left(\sum_{j=1}^m\sum_{c=1}^ns_{uj}a_{jc}t_{cv}\right)
=\left\{\begin{array}{ll}
p_i^{\alpha_{iu}} &\hbox{if}\;u=v;\\
0    &\hbox{otherwise}.
\end{array}\right.
$$
Since $\mathbb{Z}_{h}\cong\mathbb{Z}_{p_1^{s_1}}\oplus\mathbb{Z}_{p_2^{s_2}}\oplus\cdots\oplus\mathbb{Z}_{p_t^{s_t}}$,
there exists the unique element $d_u\in\mathbb{Z}_{h}$ such that $\pi_i(d_u)=p_i^{\alpha_{iu}}$ for $i=1,2,\ldots,t$,
which imply that
$$\sum_{j=1}^m\sum_{c=1}^ns_{uj}a_{jc}t_{cv}
=\left\{\begin{array}{ll}
d_u &\hbox{if}\;u=v;\\
0    &\hbox{otherwise}.
\end{array}\right.
$$
So, we have $SAT={\rm diag}\,(d_1,d_2,\ldots,d_m)$.

For each $u$ with $1\leq u\leq m$, since  $\pi_i(\prod_{i=1}^tp_i^{\alpha_{iu}})=p_i^{\alpha_{iu}}\pi_i(\prod_{j\not=i}p_j^{\alpha_{ju}})$
and $\pi_i(\prod_{j\not=i}p_j^{\alpha_{ju}})\in\mathbb{Z}_{p_i^{s_i}}^\ast$, we obtain that
$\pi_i(\prod_{i=1}^tp_i^{\alpha_{iu}})$ and $p_i^{\alpha_{iu}}=\pi_i(d_u)$ are associates for $i=1,2,\ldots,t$. By Lemma~\ref{lem2.4},
$d_u$ and $\prod_{i=1}^tp_i^{\alpha_{iu}}$ are associates for $u=1,2,\ldots,m$.
Then there exist some unit $x_u\in\mathbb{Z}_h^\ast$ such that $d_u=x_u\prod_{i=1}^tp_i^{\alpha_{iu}}$ for $u=1,2,\ldots,m$. Therefore,
$A$ is equivalent to $D$ in (\ref{equa2}). By Lemma~\ref{lem3.3}, the array $\Omega$ in (\ref{equa3})
 is uniquely determined by $A$.
$\qed$

By Theorem~\ref{lem3.4}, the vector $\Omega$  in (\ref{equa3}) is an invariance of $A$. We call it {\it invariant  factor} of $A$. 
 For a given $\Omega$ in (\ref{equa3}), let $\mathbb{Z}_{h}^{m\times n}(\Omega)$ denote the orbit of all the matrices $A\in \mathbb{Z}_{h}^{m\times n}$ with the invariant  factor $\Omega$. Let $\Omega_i=(\alpha_{i1},\alpha_{i2},\ldots,\alpha_{im})$  for $i=1,2,\ldots,t$. Then 
 $\mathbb{Z}_{p_i^{s_i}}^{m\times n}(\Omega_i)$ 
 \footnote{The length of the orbit $\mathbb{Z}_{p_i^{s_i}}^{m\times n}(\Omega_i)$  is given  in \cite{Guo2}  for $i=1,2,\ldots,t$.}
 is an orbit of $\mathbb{Z}_{p_i^{s_i}}^{m\times n}$  under the action of
 $G\!L_m(\mathbb{Z}_{p_i^{s_i}})\times G\!L_n(\mathbb{Z}_{p_i^{s_i}})$ for $i=1,2,\ldots,t$.

  \begin{thm}\label{lem3.4-1}
 Let $m\leq n$. The number of orbits of $\mathbb{Z}_{h}^{m\times n}$ under the action of $G\!L_m(\mathbb{Z}_{h})\times G\!L_n(\mathbb{Z}_{h})$ is
 $\prod_{i=1}^t{s_i+m\choose m}$.
 \end{thm}
 \proof By Theorem~\ref{lem3.4}, the number of the orbits of $\mathbb{Z}_{h}^{m\times n}$ under the action of $G\!L_m(\mathbb{Z}_{h})\times G\!L_n(\mathbb{Z}_{h})$ is equal to the number of the invariant  factors $\Omega$ in (\ref{equa3}), where  $0\leq\alpha_{i1}\leq\alpha_{i2}\leq\cdots\leq\alpha_{im}\leq s_i$ for $i=1,2,\ldots,t$.
 Note that the number of vectors $(\alpha_{i1},\alpha_{i2},\ldots,\alpha_{im})$ is equal to the number of all the $m$-multisets  from a finite set of cardinality $s_i+1$, i.e. ${s_i+m\choose m}$.
 Therefore, the number of orbits of $\mathbb{Z}_{h}^{m\times n}$ under the action of $G\!L_m(\mathbb{Z}_{h})\times G\!L_n(\mathbb{Z}_{h})$ is
 $\prod_{i=1}^t{s_i+m\choose m}$.
 $\qed$

\begin{thm}\label{lem3.4-2}
 Let $\Omega$ be as in (\ref{equa3}). Then the length of
 $\mathbb{Z}_{h}^{m\times n}(\Omega)$  is equal to
 $$\prod_{i=1}^t|\mathbb{Z}_{p_i^{s_i}}^{m\times n}(\Omega_i)|.$$
 \end{thm}
 \proof Let $\pi_i$ (resp. $\pi$) denote also the restriction of $\pi_i$ (resp. $\pi$) on $\mathbb{Z}_{h}^{m\times n}(\Omega)$ for $i=1,2,\ldots,t$. 
Let $u_{ij}=\pi_i(\prod_{c\not=i}p_c^{\alpha_{cj}})$ and $D_i={\rm diag}(u_{i1},u_{i2},\ldots,u_{im})$ for $i=1,2,\ldots,t$ and $j=1,2,\ldots,m$. By Lemma~\ref{lem3.1} and Theorem~\ref{lem3.4}, $\pi$ is a bijective map  such that
$$\begin{array}{ccc}
\pi:\mathbb{Z}_{h}^{m\times n}(\Omega) &\rightarrow &D_1\mathbb{Z}_{p_1^{s_1}}^{m\times n}(\Omega_1)\times D_2\mathbb{Z}_{p_2^{s_2}}^{m\times n}(\Omega_2)
 \times\cdots\times D_t\mathbb{Z}_{p_t^{s_t}}^{m\times n}(\Omega_t),\\
  A &\mapsto& (\pi_1(A),\pi_2(A),\ldots,\pi_t(A)).
  \end{array}$$
It follows that 
$$|\mathbb{Z}_{h}^{m\times n}(\Omega)|=\prod_{i=1}^t|D_i\mathbb{Z}_{p_i^{s_i}}^{m\times n}(\Omega_i)|=\prod_{i=1}^t|\mathbb{Z}_{p_i^{s_i}}^{m\times n}(\Omega_i)|,$$
as desired. $\qed$

\begin{lemma}\label{lem3.5}
Let $A=SDT\in \mathbb{Z}_{h}^{m\times n}$, where $S$ and $T$ are invertible, and $D$ is as in (\ref{equa2}).
Then the inner rank of $A$ is $\max\{c: (\alpha_{1c},\alpha_{2c},\ldots,\alpha_{tc})\not=(s_1,s_2,\ldots,s_t)\}$.
\end{lemma}
\proof
Let $j=\max\{c: (\alpha_{1c},\alpha_{2c},\ldots,\alpha_{tc})\not=(s_1,s_2,\ldots,s_t)\}$.
Write $S=(S_1, S_2)$ and $T=\left(\begin{array}{c}T_1\\ T_2\end{array}\right)$,
where $S_1\in \mathbb{Z}_{h}^{m\times j}$ and $T_1\in \mathbb{Z}_{h}^{j\times n}$.
Then $$A =S_1{\rm diag}\,\left(\prod_{i=1}^tp_i^{\alpha_{i1}},\prod_{i=1}^tp_i^{\alpha_{i2}},\ldots,\prod_{i=1}^tp_i^{\alpha_{ij}}\right)T_1,$$ which implies that $\rho(A)\leq j$.

Suppose $\rho(A)=\ell<j$. Then $A=A_1A_2$ with $A_1\in \mathbb{Z}_{h}^{m\times \ell}$ and $A_2\in \mathbb{Z}_{h}^{\ell\times n}$. By Theorem~\ref{lem3.4}, there exist invertible matrices $S_1,T_1,S_2,T_2$ and $\ell\times \ell$ diagonal matrices $D_1,D_2$, such that
$$A_1=S_1\left(\begin{array}{c}
D_1\\
0\end{array}\right)T_1\quad\hbox{and}\quad
A_2=S_2(D_2,0)T_2,$$
which imply that
$$A=A_1A_2=S_1\left(\begin{array}{cc}
D_1T_1S_2D_2 &\\
& 0
\end{array}\right)T_2,$$
a contradiction since $j$ is uniquely determined by $A$. Therefore, $\rho(A)=j$. $\qed$

 Let   $h_i=h/p_i^{s_i}$ and $\theta_i$  be the natural surjective homomorphism from
$\mathbb{Z}_{h}$ to $\mathbb{Z}_{h_i}$ for $i=1,2,\ldots,t$. For $A=(a_{uv})\in\mathbb{Z}_{h}^{m\times n}$, let
$\theta_i(A)=(\theta_i(a_{uv}))$ for $i=1,2,\ldots,t$.

\begin{lemma}\label{lem3.8}
Let $A\in\mathbb{Z}_{h}^{m\times n}$. Then
$$\rho(A)=\max\{\rho(\pi_i(A)): i=1,2,\ldots,t\}=\max\{\rho(\theta_i(A)): i=1,2,\ldots,t\}.$$
\end{lemma}
\proof
By Theorem~\ref{lem3.4}, there exist two invertible matrices $S$ and $T$ such that
$A=SDT,$ where $D$ is as in (\ref{equa2}). Let $\rho(A)=\ell$.
Write $S=(s_{uv})$ and $T=(t_{uv})$. For each $i$ with $1\leq i\leq t$, write $u_c=\pi_i(\prod_{j\not=i}p_j^{\alpha_{jc}})\in\mathbb{Z}_{p_i^{s_i}}^\ast$ for $c=1,2,\ldots,\ell$. Since $a_{uv}=\sum_{c=1}^{\ell}s_{uc}t_{cv}\prod_{j=1}^tp_j^{\alpha_{jc}}$, we obtain
$$
\pi_i(A)=(\pi_i(a_{uv}))
   =\pi_{i}(S)\hbox{diag}\,(p_i^{\alpha_{i1}}u_1,p_i^{\alpha_{i2}}u_2,\dots,p_i^{\alpha_{i\ell}}u_{\ell},0)\pi_{i}(T),
$$
which implies that $$\rho(\pi_i(A))=\rho(\hbox{diag}\,(p_i^{\alpha_{i1}}u_1,p_i^{\alpha_{i2}}u_2,\dots,p_i^{\alpha_{i\ell}}u_{\ell},0))\leq \ell=\rho(A).$$
Therefore, $\rho(A)\geq\max\{\rho(\pi_i(A)): i=1,2,\ldots,t\}$.
Since $(\alpha_{1\ell},\alpha_{2\ell},\ldots,\alpha_{t\ell})\not=(s_1,s_2,\ldots,s_t)$,
there exists some $j$ such that $\alpha_{j\ell}<s_j$, which implies that
$\rho(\pi_j(A))=\ell=\rho(A)$.

For each $i$ with $1\leq i\leq t$, write $v_c=\theta_i(p_i^{\alpha_{ic}})\in\mathbb{Z}_{h_i}^\ast$ for $c=1,2,\ldots,\ell$, where $h_i=h/p_i^{s_i}$. Since $a_{uv}=\sum_{c=1}^{\ell}s_{uc}t_{cv}\prod_{j=1}^tp_j^{\alpha_{jc}}$, we obtain
$$
\theta_i(A)=(\theta_i(a_{uv}))
   =\theta_{i}(S)\hbox{diag}\,\left(v_1\prod_{j\not=i}p_j^{\alpha_{j1}},v_2\prod_{j\not=i}p_j^{\alpha_{j2}},\dots,
   v_{\ell}\prod_{j\not=i}p_j^{\alpha_{j\ell}},0\right)\theta_{i}(T),
$$
which implies that $$\rho(\theta_i(A))=\rho\left(\hbox{diag}\,\left(v_1\prod_{j\not=i}p_j^{\alpha_{j1}},v_2\prod_{j\not=i}p_j^{\alpha_{j2}},\dots,
   v_{\ell}\prod_{j\not=i}p_j^{\alpha_{j\ell}},0\right)\right)\leq \ell=\rho(A).$$
Therefore, $\rho(A)\geq\max\{\rho(\theta_i(A)): i=1,2,\ldots,t\}$.
Since $(\alpha_{1\ell},\alpha_{2\ell},\ldots,\alpha_{t\ell})\not=(s_1,s_2,\ldots,s_t)$,
there exists some $j$ such that $\prod_{b\not=j}p_b^{\alpha_{b\ell}}\not=0$ in $\mathbb{Z}_{h_j}$, which implies that
$\rho(\theta_j(A))=\rho(A)$.
$\qed$

\section{Generalized bilinear forms graphs}
In this section, we always assume that $1\leq r\leq m\leq n$, and ${\rm Bil}_r$  is the generalized bilinear forms graph
${\rm Bil}_r(\mathbb{Z}_{h}^{m\times n})$. For subsets $S_1,S_2$ and $S_3$ of $\mathbb{Z}_{h}$, let
\begin{eqnarray*}
&&\left(\begin{array}{cc}
S_1^{r\times r} & S_2^{r\times(n-r)}\\
S_3^{(m-r)\times r} & 0
\end{array}\right)\\
&=&\left\{\left(\begin{array}{cc}
X_1 & X_2\\  X_3 & 0
\end{array}\right):X_1\in S_1^{r\times r},X_2\in S_2^{r\times(n-r)}, X_3\in S_3^{(m-r)\times r}\right\}.
\end{eqnarray*}

\begin{thm}\label{lem4.1}
${\rm Bil}_r$ is a connected vertex transitive graph.
\end{thm}
\proof Note that the map
$X \mapsto S^{-1}XT+A\; (X\in \mathbb{Z}_{h}^{m\times n})$
is an automorphism of ${\rm Bil}_r$ for any $S\in G\!L_{m}(\mathbb{Z}_{h}),T\in G\!L_n(\mathbb{Z}_{h})$ and $A\in\mathbb{Z}_{h}^{m\times n}$. For any fixed vertex $A$, since the bijective map $X \mapsto X+A$ is an automorphism of ${\rm Bil}_r$, ${\rm Bil}_r$ is vertex transitive.
 $\qed$

\begin{thm}{\rm(See \cite{Huang2}.)}\label{lem5.1}
The independence number, the clique number and the chromatic number
of ${\rm Bil}_r(\mathbb{Z}_{p^s}^{m\times n})$ are
$$\alpha({\rm Bil}_r(\mathbb{Z}_{p^s}^{m\times n}))=p^{sn(m-r)},\;\omega({\rm Bil}_r(\mathbb{Z}_{p^s}^{m\times n}))=p^{snr},\;
\chi({\rm Bil}_r(\mathbb{Z}_{p^s}^{m\times n}))=p^{snr}.$$
\end{thm}

If $\Gamma$ is  a vertex transitive graph, by Lemma~2.7.2 in \cite{Roberson}, we obtain the following inequality
\begin{equation}\label{equa20}
\chi(\Gamma) \geq |V(\Gamma)|/\alpha(\Gamma) \geq \omega(\Gamma).
\end{equation}

Let $G$ be a finite group and  $C$ be a subset of $G$ that is closed under taking inverses and does not contain the identity. The {\it Cayley graph} $\hbox{Cay}(G, C)$ is the graph with vertex set $G$ and two vertices $x$ and $y$ are adjacent if and only if $xy^{-1}\in C$. A Cayley graph $\hbox{Cay}(G, C)$ is {\it normal}
  if $xCx^{-1}=C$ for all $x\in G$. Clearly, ${\rm Bil}_r$ is a normal Cayley graph on the matrix additive group $G$ of $\mathbb{Z}_{h}^{m\times n}$ and the inverse closed subset $C$ is the set of all matrices of inner rank $\leq r$ excluding $0$.

  \begin{lemma}\label{lem4.11}{\rm(See Corollary~6.1.3 in \cite{Godsil}.)}
   Let $\Gamma$ be a normal Cayley graph. If $\alpha(\Gamma)\omega(\Gamma)=|V(\Gamma)|$,
    then $\chi(\Gamma)=\omega(\Gamma)$.
  \end{lemma}

\begin{thm}\label{lem5.2}
The independence number, the clique number and the chromatic number
of  ${\rm Bil}_r$ are
$$\alpha({\rm Bil}_r)=h^{n(m-r)},\;
\omega({\rm Bil}_r)=h^{nr},\;
\chi({\rm Bil}_r)=h^{nr}.$$
\end{thm}
\proof By Theorem~\ref{lem5.1}, we have
 $$\alpha({\rm Bil}_r(\mathbb{Z}_{p_i^{s_i}}^{m\times n}))=p_i^{s_in(m-r)}\;\hbox{for}\;i=1,2,\ldots,t.$$
  Write $\xi_i=p_i^{s_in(m-r)}$ for $i=1,2,\ldots,t$. Let ${\cal S}_i=\{S_{i1},S_{i2},\ldots,S_{i\xi_i}\}$ be a largest independent set of ${\rm Bil}_r(\mathbb{Z}_{p_i^{s_i}}^{m\times n})$, where $S_{i\eta_i}\in \mathbb{Z}_{p_i^{s_i}}^{m\times n}$
 for  $1\leq i\leq t$ and $1\leq \eta_i\leq \xi_i$.
 Let $${\cal S}_1\times{\cal S}_2\times\cdots\times{\cal S}_t=\{(S_{1\eta_1},S_{2\eta_2},\ldots,S_{t\eta_t}):1\leq \eta_i\leq \xi_i
 \;\hbox{for}\; i=1,2,\ldots,t\}.$$
 Then $|{\cal S}_1\times{\cal S}_2\times\cdots\times{\cal S}_t|=h^{n(m-r)}$.

 Note that there is a bijective map $\pi$ from $\mathbb{Z}_{h}^{m\times n}$ to $\mathbb{Z}_{p_1^{s_1}}^{m\times n}\times \mathbb{Z}_{p_2^{s_2}}^{m\times n}
 \times\cdots\times\mathbb{Z}_{p_t^{s_t}}^{m\times n}$ such that $\pi(A)=(\pi_1(A),\pi_2(A),\ldots,\pi_t(A))$ for every $A\in \mathbb{Z}_{h}^{m\times n}$. Let
  $${\cal S}=\{\pi^{-1}(S_{1\eta_1},S_{2\eta_2},\ldots,S_{t\eta_t}):1\leq \eta_i\leq \xi_i\;\hbox{for}\;i=1,2,\ldots,t\}.$$
  Since ${\cal S}_i$ is a largest independent set of ${\rm Bil}(\mathbb{Z}_{p_i^{s_i}}^{m\times n})$
  for $i=1,2,\ldots,t$,
  by Lemma~\ref{lem3.8}, ${\cal S}$ is a independent set of ${\rm Bil}_r$, which implies that
  $\alpha({\rm Bil}_r) \geq|{\cal S}|=h^{(m-r)n}.$
   From (\ref{equa20}), we deduce that $\omega({\rm Bil}_r)\leq h^{nr}$.
   On the other hand, $${\cal M}=\left\{\left(\begin{array}{c}
   X\\ 0
   \end{array}\right): X\in\mathbb{Z}_{h}^{r\times n}\right\}$$
   is a clique of cardinality $h^{nr}$. Therefore, we obtain
   $\omega({\rm Bil}_r)= h^{nr}$.

   From (\ref{equa20}) and $\omega({\rm Bil}_r)= h^{nr}$, we deduce that
   $$\alpha({\rm Bil}_r)
   \leq\frac{h^{mn}}{h^{nr}}=h^{n(m-r)},$$
   which implies that $\alpha({\rm Bil}_r)=h^{n(m-r)}$.
   Since ${\rm Bil}_r$ is a normal Cayley graph, by Lemma~\ref{lem4.11},
   $\chi({\rm Bil}_r)=h^{nr}$. $\qed$

   \medskip
   \noindent{\bf Remark}. As a natural extension of rank distance code over a finite field,
   an $(m\times n,r+1)$ rank distance code over $\mathbb{Z}_h$ is a subset ${\cal S}$ of $\mathbb{Z}_h^{m\times n}$
   with minimum rank distance $r+1$, i.e. $\rho(A-B)\geq r+1$ for all $A,B\in{\cal S}$ with $A\not=B$.
   Every independent set of ${\rm Bil}_r$ is an $(m\times n,r+1)$ rank distance code over $\mathbb{Z}_h$ and vice versa.
   Let ${\cal S}$ be an $(m\times n,r+1)$ rank distance code  over $\mathbb{Z}_h$. By Theorem~\ref{lem5.2}, we have $|{\cal S}|\leq h^{n(m-r)}$.
   If a rank distance code ${\cal S}$ over $\mathbb{Z}_h$ satisfies $|{\cal S}|=h^{n(m-r)}$, then ${\cal S}$ is called an
   $(m\times n,r+1)$ maximum rank distance code  over $\mathbb{Z}_h$.
   Then an $(m\times n,r+1)$ maximum rank distance code over $\mathbb{Z}_h$ is a largest independent set of ${\rm Bil}_r$ and vice versa.

   \medskip
To determine the geometric structures of maximum cliques of ${\rm Bil}_r$, we construct the following set
$$
{\cal C}_r^{m\times n}(\alpha_1,\alpha_2,\ldots,\alpha_t)=\left(\begin{array}{cc}
\mathbb{Z}_{h}^{r\times r} & J_{(\alpha_1,\alpha_2,\ldots,\alpha_t)}^{r\times(n-r)}\\
J_{(s_1-\alpha_1,s_2-\alpha_2,\ldots,s_t-\alpha_t)}^{(m-r)\times r} & 0
\end{array}\right),
$$
where either $\alpha_i=0$ or $\alpha_i=s_i$ for $i=1,2,\ldots,t$.
Then $|{\cal C}_r^{m\times n}(\alpha_1,\alpha_2,\ldots,\alpha_t)|=h^{nr}$
if $\alpha_1=\alpha_2=\cdots=\alpha_t=0$ or $m=n.$

\begin{lemma}\label{lem4.2}
Let $\alpha_1=\alpha_2=\cdots=\alpha_t=0$ or $m=n$. Then ${\cal C}_r^{m\times n}(\alpha_1,\alpha_2,\ldots,\alpha_t)$  is a maximum clique of ${\rm Bil}_r$.
\end{lemma}
\proof
By Theorem~\ref{lem5.2}, we only need to prove ${\cal C}_r^{m\times n}(\alpha_1,\alpha_2,\ldots,\alpha_t)$ is a clique of ${\rm Bil}_r$.
For any $X\in {\cal C}_r^{m\times n}(\alpha_1,\alpha_2,\ldots,\alpha_t)$, we have either
$$\pi_i(X)\in\left(\begin{array}{cc}
\mathbb{Z}_{p_i^{s_i}}^{m\times r}& 0
\end{array}\right)\quad\hbox{or}\quad
\pi_i(X)\in\left(\begin{array}{c}
\mathbb{Z}_{p_i^{s_i}}^{r\times n}\\ 0
\end{array}\right),$$
which implies that $\rho(\pi_i(X))\leq r$ for $i=1,2,\ldots,t$. By Lemma~\ref{lem3.8}, we have $\rho(X)\leq r$.
Since $X-Y\in{\cal C}_r^{m\times n}(\alpha_1,\alpha_2,\ldots,\alpha_t)$ for any $X,Y\in{\cal C}_r^{m\times n}(\alpha_1,\alpha_2,\ldots,\alpha_t)$,
${\cal C}_r^{m\times n}(\alpha_1,\alpha_2,\ldots,\alpha_t)$ is a clique of ${\rm Bil}_r$.
$\qed$

\begin{thm}{\rm(See Theorem~4.4 in \cite{Huang2}.)}\label{lem5.3}
Let ${\cal C}$ be a maximum clique of ${\rm Bil}_r(\mathbb{Z}_{p^s}^{m\times n})$.
Then ${\cal C}$ is of the form either
$$
S\left(\begin{array}{c}
\mathbb{Z}_{p^s}^{r\times n}\\ 0
\end{array}\right)+B_0,
$$
or
$$
\left(\begin{array}{cc}
\mathbb{Z}_{p^s}^{m\times r} & 0
\end{array}\right)T+B_0\;\hbox{with}\;m=n,
$$
where $S,T\in G\!L_{m}(\mathbb{Z}_{p^s})$ and $B_0\in \mathbb{Z}_{p^s}^{m\times n}$ are fixed.
\end{thm}

\begin{thm}\label{lem5.5}
Let ${\cal C}$ be a maximum clique of ${\rm Bil}_r$.
Then ${\cal C}$ is of the form either
\begin{equation}\label{equa24}
S{\cal C}_r^{m\times n}(0,0,\ldots,0)+B_0,
\end{equation}
\begin{equation}\label{equa25}
{\cal C}_r^{m\times n}(s_1,s_2,\ldots,s_t)T+B_0\;\hbox{with}\;m=n,
\end{equation}
or
\begin{equation}\label{equa26}
S{\cal C}_r^{m\times n}(\alpha_1,\alpha_2,\ldots,\alpha_t)T+B_0\;\hbox{with}\;m=n\;\hbox{and}\;(\alpha_1,\alpha_2,\ldots,\alpha_t)\not=0,(s_1,s_2,\ldots,s_t),
\end{equation}
where $S,T\in G\!L_{m}(\mathbb{Z}_{h})$ and $B_0\in \mathbb{Z}_{h}^{m\times n}$ are fixed.
\end{thm}
\proof Let ${\cal C}$ be a maximum clique of ${\rm Bil}_r(\mathbb{Z}_{h}^{m\times n})$ containing $B_0$.
If $t=1$, then this theorem holds by Theorem~\ref{lem5.3}.
From now on,  we assume that $t\geq 2$.
By Theorem~\ref{lem5.2}, we have $\omega=|{\cal C}|=h^{nr}$. By the bijection $X\mapsto X-B_0$,
 we may assume that ${\cal C}$ contains $0$.

By Lemma~\ref{lem3.8}, $\pi_i({\cal C})$ is a clique of  ${\rm Bil}_r(\mathbb{Z}^{m\times n}_{p_i^{s_i}})$.
Since $0\in{\cal C}$ and $\pi_i(0)=0$, we have $0\in\pi_i({\cal C})$.
Suppose that $\pi_i({\cal C})=\{\pi_i(A_{i1}),\pi_i(A_{i2}),\ldots,\pi_i(A_{i\omega_i})\}$. By Theorem~\ref{lem5.2},
we have $\omega_i\leq\omega({\rm Bil}_r(\mathbb{Z}_{p_i^{s_i}}^{m\times n}))=p_i^{s_inr}$.
For any vertex $X$ in $\pi_i({\cal C})$, let $\pi_i^{-1}(X)=\{Y\in{\cal C}:\pi_i(Y)=X\}$.
Then ${\cal C}$ has a partition into $\omega_i$ cliques ${\cal C}=\bigcup_{c=1}^{\omega_i}{\cal C}_c$,
where $\pi_i({\cal C}_c)=\{\pi_i(A_{ic})\}$ for $c=1,2,\ldots,\omega_i$. Therefore,
$|{\cal C}|=\sum_{c=1}^{\omega_i}|{\cal C}_c|=h^{nr}$.

Let $n_c=|{\cal C}_c|$ for $c=1,2,\ldots,\omega_i$. Then there exists a clique
$\{B_{1c},B_{2c},\ldots,B_{n_cc}\}\subseteq J_{(0,\ldots,0,s_i,0,\ldots,0)}^{m\times n}$ such that
$${\cal C}_c=\{A_{ic}+B_{1c},A_{ic}+B_{2c},\ldots,A_{ic}+B_{n_cc}\}$$
for $c=1,2,\ldots,\omega_i$.
By Lemma~\ref{lem3.8}, $\{\theta_i(B_{1c}),\theta_i(B_{2c}),\ldots,\theta_i(B_{n_cc})\}$ is a clique of
 ${\rm Bil}_r(\mathbb{Z}_{h_i}^{m\times n})$, where $h_i=h/p_i^{s_i}$.
 By Theorem~\ref{lem5.2}, we obtain
 $$n_c\leq \omega({\rm Bil}_r(\mathbb{Z}_{h_i}^{m\times n}))= h_i^{nr}
 \;\hbox{for}\;c=1,2,\ldots,\omega_i.$$
 Thus $h^{nr}=|{\cal C}|=\sum_{c=1}^{\omega_i}n_c\leq\omega_ih_i^{nr}$, which implies that
 $p_i^{s_inr}\leq\omega_i$. So, $\omega_i=p_i^{s_inr}$. Therefore, $\pi_i({\cal C})$ is a
 maximum clique of ${\rm Bil}_r(\mathbb{Z}_{p_i^{s_i}}^{m\times n})$ containing $0$.
 By Theorem~\ref{lem5.3}, we have either
 \begin{equation}\label{equa29}
\pi_i({\cal C})=S_i\left(\begin{array}{c}
\mathbb{Z}_{p_i^{s_i}}^{r\times n}\\ 0
\end{array}\right)
\quad\hbox{or}\quad
\pi_i({\cal C})=\left(\begin{array}{cc}
\mathbb{Z}_{p_i^{s_i}}^{m\times r} & 0
\end{array}\right)T_i\;\hbox{with}\;m=n,
\end{equation}
where $S_i,T_i\in G\!L_{m}(\mathbb{Z}_{p_i^{s_i}})$ are fixed.

{\bf Case}~1: $m<n$. By Lemma~\ref{lem3.1}, there exists the unique $S\in G\!L_{m}(\mathbb{Z}_{h})$ such that
$\pi_i(S)=S_i$ for $i=1,2,\ldots,t$.
For any subset ${\cal A}_i$ of $\mathbb{Z}_{p_i^{s_i}}^{m\times n}$,
let $\pi_i^{-1}({\cal A}_i)=\bigcup_{X\in{\cal A}_i}\pi_i^{-1}(X)$
for $i=1,2,\ldots,t$.
Since $\pi_i(\mathbb{Z}_{h})=\mathbb{Z}_{p_i^{s_i}}$,  we obtain
$$
\pi_i({\cal C})=\pi_i(S)\left(\begin{array}{c}
\pi_i(\mathbb{Z}_{h})^{r\times n}\\ 0
\end{array}\right)
=\pi_i\left(S\left(\begin{array}{c}
\mathbb{Z}_{h}^{r\times n}\\ 0
\end{array}\right)\right)\;\hbox{for}\;i=1,2,\ldots,t,
$$
which imply that
$${\cal C}\subseteq \pi_i^{-1}\left(\pi_i\left(S\left(\begin{array}{c}
\mathbb{Z}_{h}^{r\times n}\\ 0
\end{array}\right)\right)\right)
=S\left(\begin{array}{c}
\mathbb{Z}_{h}^{r\times n}\\ 0
\end{array}\right)$$
since $\mathbb{Z}_{h}\cong\mathbb{Z}_{p_1^{s_1}}\oplus\mathbb{Z}_{p_2^{s_2}}\oplus\cdots\oplus\mathbb{Z}_{p_t^{s_t}}$.
Since ${\cal C}$ is a maximum clique, by Lemma~\ref{lem4.2}, one obtains
$${\cal C}=S\left(\begin{array}{c}
\mathbb{Z}_{h}^{r\times n}\\ 0
\end{array}\right)
=S{\cal C}_r^{m\times n}(0,0,\ldots,0).$$
Thus the original ${\cal C}$ is of the form (\ref{equa24}).

{\bf Case}~2: $m=n$. By (\ref{equa29}),  there exists some subset $\Pi$ of $[t]=\{1,2,\ldots,t\}$ such that
$$\pi_i({\cal C})=S_i\left(\begin{array}{c}
\mathbb{Z}_{p_i^{s_i}}^{r\times n}\\ 0
\end{array}\right)\;\hbox{if}\;i\in \Pi\quad\hbox{and}\quad
\pi_i({\cal C})=\left(\begin{array}{cc}
\mathbb{Z}_{p_i^{s_i}}^{n\times r} & 0
\end{array}\right)T_i\;\hbox{if}\;i\not\in \Pi.$$
If $\Pi=[t]$, similar to the proof of Case~1, we have
${\cal C}=S{\cal C}_r^{n\times n}(0,0,\ldots,0)$, where $S\in G\!L_{n}(\mathbb{Z}_{h})$.
Thus the original ${\cal C}$ is of the form (\ref{equa24}).
If $\Pi=\emptyset$, similar to the proof of Case~1, we have
${\cal C}={\cal C}_r^{n\times n}(s_1,s_2,\ldots,s_t)T$, where $T\in G\!L_{n}(\mathbb{Z}_{h})$.
Thus the original ${\cal C}$ is of the form (\ref{equa25}).
Assume that $\Pi\not=[t],\emptyset$. By Lemma~\ref{lem3.1}, there exists the unique $(S,T)\in G\!L_{n}(\mathbb{Z}_{h})\times G\!L_{n}(\mathbb{Z}_{h})$ such that
$$\pi_i(S)=\left\{\begin{array}{ll}
S_i &\hbox{if}\;i\in \Pi;\\
I_n &\hbox{if}\;i\not\in \Pi,
\end{array}\right.
\quad\hbox{and}\quad
\pi_i(T)=\left\{\begin{array}{ll}
I_n &\hbox{if}\;i\in \Pi;\\
T_i &\hbox{if}\;i\not\in \Pi.
\end{array}\right.
$$
Without loss of generality, we may assume that $\Pi=\{1,2,\ldots,|\Pi|\}$.
Since $$\pi_i(J_{(0,\ldots,0,s_{|\Pi|+1},\ldots,s_t)})=\left\{\begin{array}{ll}
\mathbb{Z}_{p_i^{s_i}} &\hbox{if}\;i\in\Pi;\\
0  &\hbox{if}\;i\not\in\Pi,
\end{array}\right.
$$
and
$$
\pi_i(J_{(s_1,\ldots,s_{|\Pi|},0,\ldots,0)})=\left\{\begin{array}{ll}
0 &\hbox{if}\;i\in\Pi;\\
\mathbb{Z}_{p_i^{s_i}}  &\hbox{if}\;i\not\in\Pi,
\end{array}\right.
$$
  we obtain
\begin{eqnarray*}
\pi_i({\cal C})&=&\pi_i(S)\left(\begin{array}{cc}
\pi_i(\mathbb{Z}_{h})^{r\times r}  & \pi_i(J_{(0,\ldots,0,s_{|\Pi|+1},\ldots,s_t)})^{r\times(n-r)} \\
\pi_i(J_{(s_1,\ldots,s_{|\Pi|},0,\ldots,0)})^{(n-r)\times r} & 0
\end{array}\right)\pi_i(T)\\
&=&\pi_i\left(S\left(\begin{array}{cc}
\mathbb{Z}_{h}^{r\times r}  & J_{(0,\ldots,0,s_{|\Pi|+1},\ldots,s_t)}^{r\times(n-r)} \\
J_{(s_1,\ldots,s_{|\Pi|},0,\ldots,0)}^{(n-r)\times r} & 0
\end{array}\right)T\right)\;\hbox{for}\;i=1,2,\ldots,t,
\end{eqnarray*}
which imply that
\begin{eqnarray*}
{\cal C}&\subseteq& \pi_i^{-1}\left(\pi_i\left(S\left(\begin{array}{cc}
\mathbb{Z}_{h}^{r\times r}  & J_{(0,\ldots,0,s_{|\Pi|+1},\ldots,s_t)}^{r\times(n-r)} \\
J_{(s_1,\ldots,s_{|\Pi|},0,\ldots,0)}^{(n-r)\times r} & 0
\end{array}\right)T\right)\right)\\
&=&S\left(\begin{array}{cc}
\mathbb{Z}_{h}^{r\times r}  & J_{(0,\ldots,0,s_{|\Pi|+1},\ldots,s_t)}^{r\times(n-r)} \\
J_{(s_1,\ldots,s_{|\Pi|},0,\ldots,0)}^{(n-r)\times r} & 0
\end{array}\right)T
\end{eqnarray*}
since $$J_{(0,\ldots,0,s_{|\Pi|+1},\ldots,s_t)}\cong\mathbb{Z}_{p_1^{s_1}}\oplus\cdots\oplus\mathbb{Z}_{p_{|\Pi|}^{s_{|\Pi|}}}\oplus\{0\}\oplus\cdots\oplus\{0\}$$
and $$J_{(s_1,\ldots,s_{|\Pi|},0,\ldots,0)}\cong\{0\}\oplus\cdots\oplus\{0\}\oplus\mathbb{Z}_{p_{|\Pi|+1}^{s_{|\Pi|+1}}}\oplus\cdots\oplus\mathbb{Z}_{p_{t}^{s_t}}.$$
Since ${\cal C}$ is a maximum clique, by Lemma~\ref{lem4.2}, one obtains
$${\cal C}=S\left(\begin{array}{cc}
\mathbb{Z}_{h}^{r\times r}  & J_{(0,\ldots,0,s_{|\Pi|+1},\ldots,s_t)}^{r\times(n-r)} \\
J_{(s_1,\ldots,s_{|\Pi|},0,\ldots,0)}^{(n-r)\times r} & 0
\end{array}\right)T
=S{\cal C}_r^{n\times n}(0,\ldots,0,s_{|\Pi|+1},\ldots,s_t)T.$$
Thus the original ${\cal C}$ is of the form (\ref{equa26}). $\qed$

By Theorem~\ref{lem5.5}, we obtain the Erd\H{o}s-Ko-Rado theorem for $\mathbb{Z}_h^{m\times n}$.

\begin{thm}\label{lem5.6}
Let ${\cal F}$ be an $r$-intersecting family of $\mathbb{Z}_h^{m\times n}$.
Then $|{\cal F}|\leq h^{nr},$
and equality holds if and only if either {\rm (a)}
${\cal F}=S{\cal C}_r^{m\times n}(0,0,\ldots,0)+B_0,$
{\rm (b)} ${\cal F}={\cal C}_r^{m\times n}(s_1,s_2,\ldots,s_t)T+B_0$ with $m=n$,
or {\rm (c)}
${\cal F}=S{\cal C}_r^{n\times n}(\alpha_1,\alpha_2,\ldots,\alpha_t)T+B_0$ with
$m=n$ and $(\alpha_1,\alpha_2,\ldots,\alpha_t)\not=0,(s_1,s_2,\ldots,s_t),$
where $S,T\in G\!L_{m}(\mathbb{Z}_{h})$ and $B_0\in \mathbb{Z}_{h}^{m\times n}$ are fixed.
\end{thm}

Now, we discuss cores of both $\hbox{Bil}_r(\mathbb{Z}_{h}^{m\times n})$  and its complement.

\begin{lemma}\label{lem5.7}{\rm(See Lemma~2.5.9 in \cite{Roberson}.)}
Let $\Gamma$ be a graph. Then the core of $\Gamma$ is the complete graph $K_c$ if and only if $\omega(\Gamma)=c=\chi(\Gamma)$.
\end{lemma}

\begin{cor}\label{lem5.8}
The core of ${\rm Bil}_r$ is a maximum clique, and ${\rm Bil}_r$ is a core if and only if $r=m$.
\end{cor}
\proof
By Theorem~\ref{lem5.2}, we have $\chi({\rm Bil}_r)=\omega({\rm Bil}_r)=h^{nr}$, which implies that
 the core of ${\rm Bil}_r$ is a maximum clique by Lemma~\ref{lem5.7}.
 Since ${\rm Bil}_r$ is a clique if and only if $r=m$, ${\rm Bil}_r$ is a core if and only if $r=m$. $\qed$

\begin{thm}\label{lem5.9}
Let $\Sigma_r$ denote the complement of ${\rm Bil}_r$.
Then $\chi(\Sigma_r)=\alpha(\Gamma)=\omega(\Sigma_r)$. Moreover,
the core of $\Sigma_r$ is a maximum clique, and $\Sigma_r$ is not a core.
\end{thm}
\proof
Note that $\alpha({\rm Bil}_r)=\omega(\Sigma_r),\alpha(\Sigma_r)=\omega({\rm Bil}_r)$,
$\chi(\Sigma_r)\geq\omega(\Sigma_r)=\alpha({\rm Bil}_r)$.
Let ${\cal S}=\{S_1,S_2,\ldots,S_\alpha\}$ be a largest independent set of ${\rm Bil}_r$,
where $\alpha=h^{(m-r)n}$ and $S_i\in\mathbb{Z}_h^{m\times n}$ for $i=1,2,\ldots,\alpha$.
Let ${\cal C}_i=S_i+{\cal C}_r^{m\times n}(0,0,\ldots,0)$ for $i=1,2,\ldots,\alpha$.
Then ${\cal C}_1,{\cal C}_2,\ldots,{\cal C}_\alpha$ are $\alpha$ maximum cliques of ${\rm Bil}_r$,
and ${\cal C}_i\cap{\cal C}_j=\emptyset$ for all $i\not=j$.
By Theorem~\ref{lem5.2}, we have $|V({\rm Bil}_r)|=\alpha\omega({\rm Bil}_r)$. So, $V({\rm Bil}_r)$ has a partition into $\alpha$ maximum cliques:
$V({\rm Bil}_r)=\bigcup_{i=1}^\alpha{\cal C}_i$. Since ${\cal C}_i$ is a largest independent set of $\Sigma_r$ for $i=1,2,\ldots,\alpha$,
 $V(\Sigma_r)$ has a partition into $\alpha$ largest independent sets: $V(\Sigma_r)=\bigcup_{i=1}^\alpha{\cal C}_i$.
  It follows that $\chi(\Sigma_r)\leq\alpha$, which implies that $\chi(\Sigma_r)=\alpha$ by $\chi(\Sigma_r)\geq\alpha$.
Similar to the proof of Corollary~\ref{lem5.8}, the core of $\Sigma_r$ is a maximum clique, and $\Sigma_r$ is not a core
 since $\Sigma_r$ is not a clique. $\qed$

\section*{Acknowledgment}
This research is supported by National Natural Science Foundation of China (11971146).

\end{document}